# Global Stabilization of a Class of Nonlinear Reaction-Diffusion PDEs by Boundary Feedback

**Iasson Karafyllis[*] and Miroslav Krstic[**]**

[*]Dept. of Mathematics, National Technical University of Athens, Zografou Campus, 15780, Athens, Greece, email: iasonkar@central.ntua.gr

[**]Dept. of Mechanical and Aerospace Eng., University of California, San Diego, La Jolla, CA 92093-0411, U.S.A., email: krstic@ucsd.edu

**Abstract**

This paper provides global exponential stabilization results by means of boundary feedback control for 1-D nonlinear unstable reaction-diffusion Partial Differential Equations (PDEs) with nonlinearities of superlinear growth. The class of systems studied are parabolic PDEs with nonlinear reaction terms that provide "damping" when the norm of the state is large (the class includes reaction-diffusion PDEs with polynomial nonlinearities). The case of Dirichlet actuation at one end of the domain is considered and a Control-Lypunov Functional construction is applied in conjunction with Stampacchia's truncation method. The paper also provides several important auxiliary results; among which is an extension of Wirtinger's inequality, used here for the construction of the Control Lyapunov functional.

## 1. Introduction

The stabilization of equilibrium profiles for parabolic Partial Differential Equations (PDEs) by means of feedback control is an important problem in control theory which has attracted the interest of many researchers. Research has focused on linear parabolic PDEs, where several methodologies are available; see [6,11] and the backstepping design methodology described in [34]. Unstable nonlinear PDEs are usually treated by performing a feedback design for the linearized model (see for example [36]). Very few feedback design methodologies have been proposed for unstable nonlinear parabolic PDEs: see the extension of the backstepping boundary feedback design in [37,38] as well as feedback designs for distributed inputs in [9,16,31,32]. In many cases the stabilization results are local, guaranteeing exponential stability in specific spatial norms. For nonlinear hyperbolic PDEs, feedback stabilization results are given in [5,10,13,14,21,30].

The recent paper [22] provided global exponential stabilization results in the $L^2$ spatial norm for 1-D nonlinear parabolic PDEs by means of boundary feedback. A small-gain methodology (based on the results in [20]) was applied to parabolic PDEs with nonlinearities that satisfy a linear growth condition. It is well-known that global stabilization by means of boundary feedback may not be possible for parabolic PDEs with nonlinearities of superlinear growth due to the loss of controllability and the existence of finite escape times (see [12,15]).

The present work is the first paper that provides global exponential stabilization results by means of boundary feedback control for 1-D nonlinear reaction-diffusion PDEs with nonlinearities of superlinear growth. The class of systems studied is a specific class of PDEs where the nonlinear terms provide "damping" when the norm of the state is large (which includes generalized Chafee-Infante PDEs and reaction-diffusion PDEs with polynomial nonlinearities studied in [1,8,33,35]). We consider the case of Dirichlet actuation at one end of the domain. Our methodology is different from the small-gain methodology in [22]: here we construct Control Lyapunov Functionals (CLFs) which provide exponential stability estimates in the $L^2$ spatial norm. Lyapunov functionals for reaction-diffusion PDEs were recently proposed in [27,28]. However, due to the fact that the interval of existence of the solution for reaction-diffusion PDEs is related to the sup-norm of the state, a more involved analysis is required in order to guarantee existence of solutions for all times. Applying Stampacchia's truncation method (as presented in the proof of Theorem 10.3 in [7]) and exploiting the "damping" of the nonlinear reaction term when the norm of the state is large, we are in a position to get bounds of the spatial sup-norm. From a control-theoretic point of view this feature indicates a sharp contrast with the nonlinear finite-dimensional case: in the nonlinear infinite-dimensional case a CLF feedback design may not be sufficient for establishing existence of solutions for all times (global solutions) and consequently may not allow a valid (not merely formal) derivation of stability properties. Additional analysis may be required in order to obtain bounds that guarantee the existence of the solution for all times and the pointwise convergence of the solution to the desired equilibrium point (important in practice).

The obtained results guarantee global exponential stability in the $L^2$ spatial norm, global exponential convergence in the spatial sup-norm and Lagrange stability as well as Lyapunov (global) stability in the $H^1$ spatial norm. Thus we obtain different stabilization results in different spatial norms of the state. This feature is well-known in the dynamical systems literature of PDEs (see for example the different attractors in different spaces presented in [35]).

We consider this work to be somewhat of a breakthrough in global stabilization of open-loop unstable equilibria for nonlinear parabolic PDEs with nonlinearities of superlinear growth. The paper advances a string of efforts over the past two decades on this topic, which have so far resulted only in either global asymptotic stabilization of equilibria that are already open-loop neutrally stable for sufficiently high viscosity-like coefficients, such as for the viscous Burgers equation [2,23,26], KdV equation [3], Kuramoto-Sivashinsky equation [25], and 2D Navier-Stokes channel flow [4], or in merely regional stabilization for the open-loop unstable viscous Burgers equation in [24] and for a class of integro-differential parabolic systems with spatial Volterra operators on the right hand side [37,38]. In this paper we make a step forward in both allowing open-loop instability and in achieving global asymptotic stabilization in closed loop.

The structure of the paper is as follows. Section 2 is devoted to the presentation of the problem and the statement of the main result (Theorem 2.1). A discussion of the main result is also provided in Section 2 as well as an example (Example 2.4), which illustrates how easily the main result can be applied to a given reaction-diffusion PDE. Section 3 of the present work is devoted to the presentation of fundamental existence-uniqueness results for semilinear reaction-diffusion PDEs with possible non-local terms. Specific results are provided for the heat equation with Dirichlet boundary conditions (Theorem 3.5). The existence-uniqueness results of Section 3 are slightly different from similar results in [8] and can be used for control purposes (since the boundary feedback gives a non-local term in the PDE after homogenizing the boundary conditions). The proof of the main result is given in Section 4 of the paper, where additional auxiliary results are also presented. Among the auxiliary results, there is an extension (Proposition 4.1) of the well-known Wirtinger's inequality (see [17]). The concluding remarks of the paper are provided in Section 5.

**Notation.** Throughout this paper, we adopt the following notation.

* $\Re_+ := [0, +\infty)$.
* Let $A \subseteq \Re^n$ be an open set and let $\Omega \subseteq \Re$ and $A \subseteq U \subseteq \bar{A}$ be given sets. By $C^0(U)$ (or $C^0(U; \Omega)$), we denote the class of continuous mappings on $U$ (which take values in $\Omega$). By $C^k(U)$ (or $C^k(U; \Omega)$), where $k \geq 1$, we denote the class of continuous functions on $U$, which have continuous derivatives of order $k$ on $U$ (and also take values in $\Omega$). For a differentiable function $u \in C^0([0,1])$, $u'(x)$ denotes the derivative with respect to $x \in [0,1]$.

* For $u,v \in L^2(0,1)$, $\langle u,v \rangle$ denotes the inner product in $L^2(0,1)$. For $p \in [1,+\infty]$ and $u \in L^p(0,1)$, $\|u\|_p$ is the norm of $L^p(0,1)$.

* Let $u: \Re_+ \times [0,1] \to \Re$ be given. We use the notation $u[t]$ to denote the profile at certain $t \geq 0$, i.e., $(u[t])(x) = u(t,x)$ for all $x \in [0,1]$. When $u(t,x)$ is (twice) differentiable with respect to $x \in [0,1]$, we use the notation $u_x(t,x)$ ($u_{xx}(t,x)$) for the (second) derivative of $u$ with respect to $x \in [0,1]$, i.e., $u_x(t,x) = \frac{\partial u}{\partial x}(t,x)$ ($u_{xx}(t,x) = \frac{\partial^2 u}{\partial x^2}(t,x)$). When $u(t,x)$ is differentiable with respect to $t$, we use the notation $u_t(t,x)$ for the derivative of $u$ with respect to $t$, i.e., $u_t(t,x) = \frac{\partial u}{\partial t}(t,x)$.

* For an integer $k \geq 1$, $H^k(0,1)$ ($H_0^k(0,1)$) denotes the Sobolev space of functions in $L^2(0,1)$ with all its weak derivatives up to order $k \geq 1$ in $L^2(0,1)$ (with $u(0) = u(1) = 0$).

## 2. Problem Statement and Main Result

Consider the scalar semilinear 1-D reaction-diffusion PDE

$$u_t = p u_{xx} + F(u), \quad (2.1)$$

where $x \in (0,1)$, $u(t,x) \in \Re$ and $p > 0$ is a constant, with Dirichlet boundary conditions

$$\begin{aligned} u(t,0) &= 0 \\ u(t,1) &= U(t) \end{aligned} \quad (2.2)$$

where $U(t) \in \Re$ is the boundary control input. The mapping $F: C^0([0,1]) \to C^0([0,1])$ that contains the reaction terms is assumed to satisfy the following conditions:

$$(F(u))(x) = f(x, u(x)), \text{ for } x \in [0,1] \text{ and } u \in C^0([0,1]) \quad (2.3)$$

where $f \in C^0([0,1] \times \Re)$ is a locally Lipschitz function that satisfies

$$u(f(x,u) - qu) \leq \gamma u^2 - B|u|^b, \text{ for all } u \in \Re, \ x \in [0,1] \quad (2.4)$$

$$|f(x,u) - qu| \leq \gamma |u| + \delta |u|^{b-1}, \text{ for all } u \in \Re \quad (2.5)$$

for certain real constants $q \in \Re$, $\gamma, \delta, B \geq 0$ and $b \geq 2$. This is the class of reaction-diffusion PDEs, for which (as implied by (2.4)) the nonlinear reaction term plays a "stabilizing role" when the state becomes "large". This class includes the polynomial reaction-diffusion PDEs studied in [1,8,33,35].

We consider the problem of existence/design of a linear boundary feedback law of the form

$$U(t) = -r \langle k, u[t] \rangle \quad (2.6)$$

where $k \in C^2([0,1])$ is a function and $r \in \Re$ is a constant, that achieves global exponential stabilization of system (2.1), (2.2) in the $L^2$ norm as well as global exponential convergence of the state to 0 in the $L^\infty$ norm.

The reasons that force us to consider only a linear (and not a nonlinear) boundary feedback law are: **(i)** for a nonlinear boundary feedback law additional technical issues arise which make the analysis of the closed-loop system almost impossible; these technical issues are explained in detail in what follows (see Remark 3.6 and Remark 4.6 below), and **(ii)** a linear boundary feedback law has the nice feature that when the nonlinear term is of the form $f(x,u) = qu - B\,\text{sgn}(u)|u|^{b-1}$, where $q \in \Re$, $B > 0$, $b > 2$ are constants, then a simple scaling argument shows that stabilization is achieved for all values of $B > 0$ and therefore, we do not need to know the exact value of the constant $B > 0$ (robust stabilization).

The statement of the main result of the paper follows.

**Theorem 2.1:** *Consider the closed-loop system (2.1), (2.2), (2.6), where $k \in C^2([0,1])$ is a function with $k(0) = 0$, $k(1) = 1$ that satisfies $k''(x) = \mu k(x)$ for $x \in [0,1]$ and for some constant $\mu \in \Re$, $p > 0$, $r \in \Re$ are constants with $1 + \|k\|_2^2 r > 0$, and $F : C^0([0,1]) \to C^0([0,1])$ is defined by (2.3) for a locally Lipschitz function $f \in C^0([0,1] \times \Re)$ that satisfies (2.4), (2.5) for certain real constants $q \in \Re$, $\gamma, \delta, B \geq 0$ and $b \geq 2$ with $B \geq \delta |r| \|k\|_b \|k\|_{\frac{b}{b-1}}$. Suppose that there exists a constant $\varepsilon \geq 0$ that satisfies the following inequality*

$$\max\left(0, r^2\left(\|k'\|_2^2 + (3\varepsilon - 1)\pi^2\|k\|_2^2 - k'(1)\right) + r\left(2(3\varepsilon - 1)\pi^2 + \frac{q}{p} - \mu\right) + \frac{3\pi^2 \varepsilon (1+\varepsilon)}{(1+\varepsilon)\|k\|_2^2 - k_1^2}\right)$$
$$< \frac{(1+\varepsilon)\|k\|_2^2 + (3\varepsilon - 1)k_1^2}{\left((1+\varepsilon)\|k\|_2^2 - k_1^2\right)\|k\|_2^2}\pi^2 - \frac{q + \gamma\left(1 + |r|\|k\|_2^2\right)}{p\|k\|_2^2} \tag{2.7}$$

*where $k_1 = \sqrt{2} \int_0^1 k(x)\sin(\pi x)dx$. Moreover, if $q > 0$ then suppose that $B > 0$ and $b > 2$. Then for every $u_0 \in H^2(0,1)$ with $u_0(0) = 0$, $u_0(1) = -r\langle k, u_0\rangle$, there exists a unique mapping $u \in C^0(\Re_+ \times [0,1]) \cap C^1((0,+\infty); L^2(0,1))$ with $u[0] = u_0$, $u[t] \in H_0^2(0,1)$ for all $t \geq 0$ for which equation (2.1) holds for all $t > 0$ and equation (2.6) holds for all $t \geq 0$. Moreover, there exist constants $G, \sigma > 0$ and a non-decreasing function $\psi : \Re_+ \to \Re_+$ such that the following estimates*

$$\|u[t]\|_2 \leq G \exp(-\sigma t)\|u[0]\|_2 \tag{2.8}$$

$$\|u[t]\|_\infty \leq \sqrt{2G}\exp(-\sigma t/2)\sqrt{\|u_0\|_2 \left(\|u_0'\|_2 + \psi(M)\|u_0\|_2\right)} \tag{2.9}$$

$$\|u_x[t]\|_2 \leq \|u_0'\|_2 + \psi(M)\|u_0\|_2 \tag{2.10}$$

*hold for all $t \geq 0$, where $M := \max\{\bar{K}, \|u_0\|_\infty, G|r|\|k\|_2 \|u_0\|_2\}$, $\bar{K} := \left(B^{-1}q\right)^{1/(b-2)}$ when $q > 0$ and $\bar{K} := 0$ when $q \leq 0$.*

**Remark 2.2: (a)** The state space for the closed-loop system (2.1), (2.2), (2.6) is the space $X = \{u \in H^2(0,1) : u(0) = u(1) + r\langle k, u\rangle = 0\}$. Notice that if $u[0] = u_0 \in X$ then $u[t] \in X$ for all $t \geq 0$.
**(b)** Estimate (2.8) shows Global Exponential Stability in the $L^2$ norm. Estimate (2.10) in conjunction with estimate (2.8) and the fact that $\|u\|_\infty \leq \|u'\|_2$, for all $u \in H^1(0,1)$ with $u(0) = 0$ allow the derivation of the estimate

$$\|u_x[t]\|_2 + \|u[t]\|_2 \leq \left(\|u_0'\|_2 + \|u_0\|_2\right)\left(1+G+\psi\left(\bar{K}+\left(1+G|r|\|k\|_2\right)\left(\|u_0'\|_2 + \|u_0\|_2\right)\right)\right), \text{ for } t \geq 0 \quad (2.11)$$

The above estimate shows that the equilibrium point $0 \in X$ is Lagrange and Lyapunov stable in the $H^1$ norm.

**(c)** Estimate (2.9) does not allow us to conclude Global Asymptotic Stability in the $L^\infty$ norm because the right hand side of (2.9) cannot be bounded by a function of the $L^\infty$ norm of the initial condition. On the other hand, estimate (2.9) is useful because it demonstrates exponential convergence in the $L^\infty$ norm. Moreover, estimate (2.9) shows Global Asymptotic Output Stability (see [18]) for the output map $X \ni u \to y = u \in L^\infty(0,1)$, which simply maps the state to a different space.

**(d)** Due to the fact that $f \in C^0([0,1] \times \Re)$ is a locally Lipschitz function, it follows that there exists a non-decreasing function $L: \Re_+ \to \Re_+$ for the mapping $F: C^0([0,1]) \to C^0([0,1])$ defined by (2.3), for which the following property holds:
$$\|F(u)-F(v)\|_\infty \leq L\left(\max\left(\|u\|_\infty, \|v\|_\infty\right)\right)\|u-v\|_\infty, \text{ for all } u,v \in C^0([0,1])$$
The above Lipschitz property will be used heavily for the existence/uniqueness results of the following section.

**Remark 2.3:** **(i)** The proof of Theorem 2.1 is based on the use of a Control Lyapunov Functional (CLF) for system (2.1), (2.2), (2.6), namely the quadratic functional

$$V(u) := \frac{1}{2}\|u\|_2^2 + \frac{r}{2}\langle k,u\rangle^2, \text{ for } u \in L^2(0,1) \quad (2.12)$$

**(ii)** It is interesting to note that for the case $r=0$ (open-loop system), inequality (2.7) for $\varepsilon=0$ gives the inequality $q < p\pi^2$, which is exactly the condition that guarantees global exponential stability in the $L^2$ norm for the linearization of system (2.1), (2.2), (2.6) around zero.

**(iii)** The linearization of system (2.1), (2.2), (2.6) around zero is only "mildly unstable". To understand this notice that $\dfrac{(1+\varepsilon)\|k\|^2 + (3\varepsilon-1)k_1^2}{(1+\varepsilon)\|k\|^2 - k_1^2} < 4$ for every non-zero $k \in L^2(0,1)$ and for every $\varepsilon \geq 0$. Therefore, inequality (2.7) (in conjunction with the fact that $\gamma \geq 0$) implies the following inequality

$$q < 4\pi^2 p$$

The above inequality implies that *at most* one mode is unstable for the linearization of system (2.1), (2.2), (2.6) around zero.

**(iv)** The functions $k \in C^2([0,1])$ with $k(0)=0$, $k(1)=1$ that satisfy $k''(x) = \mu k(x)$ for $x \in [0,1]$ and for some constant $\mu \in \Re$ are:

- $k(x) = x$, where $\mu = 0$,
- $k(x) = \dfrac{\sinh(cx)}{\sinh(c)}$, where $\mu = c^2$, $c > 0$,
- $k(x) = \dfrac{\sin(\omega x)}{\sin(\omega)}$, where $\mu = -\omega^2$, $\omega > 0$ and $\omega \neq n\pi$ for $n = 1, 2, ...$

Therefore, one of the above functions can be used as the kernel of the boundary feedback stabilizer $u(t,1) = -r\langle k, u[t]\rangle$.

The following example illustrates how easily we can use Theorem 2.1 for the global exponential stabilization of a reaction-diffusion PDE.

**Example 2.4:** Consider the nonlinear reaction-diffusion PDE:

$$\begin{aligned} u_t &= pu_{xx} + qu - Bu^3 \\ u(t,0) &= 0 \\ u(t,1) &= U(t) \end{aligned} \quad (2.13)$$

for $(t,x) \in (0,+\infty) \times (0,1)$, where $p > 0$, $q, B \in \Re$ are constants. We show next that global stabilization of system (2.13) can be achieved provided that

$$B > 0, \quad \frac{q}{p\pi^2} < (1-g)(1+s) \approx 1.34 \quad (2.14)$$

where $s = \dfrac{7\pi^2 - 18 - \sqrt{(7\pi^2 - 18)^2 - 72\pi^2}}{2\pi^2} \approx 0.38$, $r = 5^{1/4}\left(\dfrac{7}{3}\right)^{3/4} \approx 2.823$, $g = \dfrac{3-r}{6} \approx 0.03$, by means of the boundary feedback law

$$U(t) = -r \int_0^1 x u(t,x) dx \quad (2.15)$$

Therefore, we are in a position to achieve global stabilization without knowing the values of the constants $B, q, p$. It should be noticed that inequalities (2.14) allow the case where $q > p\pi^2$, i.e. the case where the linearization of system (2.13) around zero with $U \equiv 0$ (open-loop system) is unstable.

Indeed, notice that the closed-loop system (2.13), (2.15) is of the form (2.1), (2.2), (2.6) with $f(x,u) = qu - Bu^3$ and $k(x) = x$ for $x \in [0,1]$. Inequalities (2.4), (2.5) hold with $\gamma = 0$, $b = 4$ and $\delta = B$. Condition (2.7) takes the form

$$\frac{3q}{p\pi^2} + \max\left(0, (3\varepsilon - 1)\frac{r^2}{3} + 2(3\varepsilon - 1)r + \frac{q}{p\pi^2}r + \frac{9\pi^2\varepsilon(1+\varepsilon)}{(1+\varepsilon)\pi^2 - 6}\right) < 3\frac{(1+\varepsilon)\pi^2 + 6(3\varepsilon - 1)}{(1+\varepsilon)\pi^2 - 6}$$

Setting $\varepsilon = \dfrac{1-s}{6}$ and using the fact that $\dfrac{18(1-s)}{(7-s)\pi^2 - 36} = s$, the above condition takes the form:

$$\frac{q}{p\pi^2(1+s)} + \frac{1}{6}\max\left(0, -\frac{r^2}{3} - 2\frac{p\pi^2(1+s) - q}{p\pi^2(1+s)}r + 3\right) < 1$$

The fact that $\dfrac{1}{6}\max\left(0, -\dfrac{r^2}{3} - 2gr + 3\right) \le g$ in conjunction with (2.14) guarantee that the above inequality (and equivalently condition (2.7)) holds. Condition $1 + \|k\|_2^2 r > 0$ holds automatically. Finally, the condition $B \ge \delta |r| \|k\|_b \|k\|_{\frac{b}{b-1}}$ with $b = 4$ and $\delta = B$ is equivalent to the condition $r \le 5^{1/4}\left(\dfrac{7}{3}\right)^{3/4}$, which also holds. Therefore, all assumptions of Theorem 2.1 hold and consequently, there exist constants $G, \sigma > 0$ and a non-decreasing function $\psi : \Re_+ \to \Re_+$ such that estimates (2.8), (2.9), (2.10) hold for all $t \ge 0$. ◁

## 3. Existence and a Uniqueness Result

The goal of this section is the proof of the following existence/uniqueness result.

**Theorem 3.1:** *Let $p > 0$ be a constant and let $\bar{F}: C^0([0,1]) \to C^0([0,1])$ be a continuous mapping with $\bar{F}(0) = 0$, for which there exists a non-decreasing function $L: \Re_+ \to (0, +\infty)$ such that the following property holds:*

$$\left\| \bar{F}(u) - \bar{F}(v) \right\|_\infty \leq L\left( \max\left( \|u\|_\infty, \|v\|_\infty \right) \right) \|u - v\|_\infty, \text{ for all } u, v \in C^0([0,1]) \tag{3.1}$$

*Let $\bar{k} \in C^2([0,1])$ be a function with $\bar{k}(0) = 0$. Then for every $u_0 \in H^2(0,1)$ with $u_0(0) = 0$, $u_0(1) = -\langle \bar{k}, u_0 \rangle$, there exists $t_{\max}(u_0) \in (0, +\infty]$ and a unique mapping $u \in C^0([0, t_{\max}(u_0)) \times [0,1]) \cap C^1((0, t_{\max}(u_0)); L^2(0,1))$ with $u[0] = u_0$, $u[t] \in H^2(0,1)$ for all $t \in [0, t_{\max}(u_0))$, for which the mapping $t \to \|u_x[t] - \bar{k}u[t]\|_2^2$ is $C^1$ on $[0, t_{\max}(u_0))$ and for which the following equations hold:*

$$u_t[t] = pu_{xx}[t] + \bar{F}(u[t]), \text{ for } t \in (0, t_{\max}(u_0)) \tag{3.2}$$

$$\begin{aligned} u(t,0) &= 0 \\ u(t,1) &= \langle \bar{k}, u[t] \rangle \end{aligned}, \text{ for } t \in [0, t_{\max}(u_0)) \tag{3.3}$$

$$\frac{d}{dt} \|u_x[t] - \bar{k}u[t]\|_2^2 = -2p \|u_{xx}[t] - \bar{k}u_x[t] - \bar{k}'u[t]\|_2^2 - 2\langle K\bar{F}(u[t]) + pGu[t], u_{xx}[t] - \bar{k}u_x[t] - \bar{k}'u[t] \rangle,$$

$$\text{for } t \in [0, t_{\max}(w_0)) \tag{3.4}$$

*where $K, G: L^2(0,1) \to L^2(0,1)$ are the continuous linear operators defined by the following equations for all $u \in L^2(0,1)$, $x \in [0,1]$:*

$$(Ku)(x) = u(x) - \int_0^x \bar{k}(s)u(s)\,ds \tag{3.5}$$

$$(Gu)(x) = 2\bar{k}'(x)u(x) - \int_0^x \bar{k}''(s)u(s)\,ds \tag{3.6}$$

*Moreover, if $t_{\max}(u_0) < +\infty$, then $\lim_{t \to t_{\max}^-} \left( \|u[t]\|_\infty \right) = +\infty$.*

In order to prove Theorem 3.1, we need to show first several auxiliary results. We follow the methodology of Chapter 4 in [8] and we start with the following lemma.

**Lemma 3.2:** *Let $p > 0$ be a constant and let $\bar{F}: C^0([0,1]) \to C^0([0,1])$ be a continuous mapping with $\bar{F}(0) = 0$, for which there exists a non-decreasing function $L: \Re_+ \to (0, +\infty)$ such that (3.1) holds. Then for every $w_0 \in H_0^2(0,1)$, $M > 0$ and for every $T > 0$ with*

$$T \leq \frac{9p^2\pi^4}{4} \left( \sqrt{2}L\left(\sqrt{2}\|w_0\|_\infty + M\right) \sum_{n=1}^\infty n^{-4/3} \left( \max\left(2, \left(\sqrt{2}\|w_0\|_\infty + M\right)M^{-1/3}\right) \right) \right)^{-3}$$

*there exists a unique mapping $w \in S$, where $S := \left\{ u \in C^0([0,T] \times [0,1]) : \max_{0 \leq t \leq T} \left( \|u[t]\|_\infty \right) \leq \sqrt{2}\|w_0\|_\infty + M \right\}$ such that the following equation holds for all $(t,x) \in [0,T] \times [0,1]$:*

$$w(t,x) = 2\sum_{n=1}^{\infty}\left[\exp\left(-pn^2\pi^2 t\right)\int_0^1 w_0(s)\sin(n\pi s)ds\right]\sin(n\pi x)$$
$$+2\sum_{n=1}^{\infty}\left[\int_0^t \exp\left(-pn^2\pi^2(t-\tau)\right)\left(\int_0^1 (\bar{F}(w[\tau]))(s)\sin(n\pi s)ds\right)d\tau\right]\sin(n\pi x) \qquad (3.7)$$

**Proof:** Let $w_0 \in H_0^2(0,1)$, $M > 0$, $T > 0$ with
$T \leq \frac{9p^2\pi^4}{4}\left(\sqrt{2}L\left(\sqrt{2}\|w_0\|_\infty + M\right)\sum_{n=1}^{\infty} n^{-4/3}\left(\max\left(2,\left(\sqrt{2}\|w_0\|_\infty + M\right)M^{-1/3}\right)\right)\right)^{-3}$ be given. Notice that

$$\bar{w}(t,x) = 2\sum_{n=1}^{\infty}\left[\exp\left(-pn^2\pi^2 t\right)\int_0^1 w_0(s)\sin(n\pi s)ds\right]\sin(n\pi x), \text{ for } (t,x)\in[0,T]\times[0,1] \qquad (3.8)$$

is the solution $\bar{w}\in C^0([0,T]\times[0,1])\cap C^1((0,T]\times[0,1])$ with $\bar{w}[t]\in C^2([0,1])$ for $t>0$ of the heat equation $\bar{w}_t[t] = p\bar{w}_{xx}[t]$ with Dirichlet boundary conditions and initial condition $\bar{w}[0] = w_0$, which by virtue of Remark 2.6 in [19,20] satisfies

$$\max_{0\leq t\leq T}\left(\|\bar{w}[t]\|_\infty\right) \leq \sqrt{2}\|w_0\|_\infty \qquad (3.9)$$

Indeed, (3.9) is established by a density argument since Remark 2.6 applies to initial conditions in $C^2([0,1])$ while here we have $w_0 \in H^2(0,1)$. For each fixed $w_0 \in H^2(0,1)$ define the operator $P:C^0([0,T]\times[0,1]) \to C^0([0,T]\times[0,1])$ by means of the equation

$$(Pw)(t,x) = \bar{w}(t,x) + 2\sum_{n=1}^{\infty}\left[\int_0^t \exp\left(-pn^2\pi^2(t-\tau)\right)\left(\int_0^1 (\bar{F}(w[\tau]))(s)\sin(n\pi s)ds\right)d\tau\right]\sin(n\pi x),$$
$$\text{for } (t,x)\in[0,T]\times[0,1] \qquad (3.10)$$

Indeed, $P:C^0([0,T]\times[0,1])\to C^0([0,T]\times[0,1])$ is a well-defined mapping by (3.10) since the facts that $2\int_0^1 \sin^2(n\pi x)dx = 1$ for $n=1,2,\ldots$, $F(0)=0$ and (3.1) imply the following inequality:

$$\left|2\int_0^t \exp\left(-pn^2\pi^2(t-\tau)\right)\left(\int_0^1 (\bar{F}(w[\tau]))(s)\sin(n\pi s)ds\right)d\tau\right| \leq \sqrt{2}\frac{L\left(\max_{0\leq t\leq T}\left(\|w[t]\|\right)\right)}{pn^2\pi^2}\max_{0\leq t\leq T}\left(\|w[t]\|\right),$$
$$\text{for all } t\in[0,T] \text{ and } n=1,2,\ldots \qquad (3.11)$$

which shows that the series in the right hand side of (3.10) converges absolutely and uniformly. Defining $S := \left\{w\in C^0([0,T]\times[0,1]): \max_{0\leq t\leq T}\left(\|w[t]\|_\infty\right)\leq \sqrt{2}\|w_0\|_\infty + M\right\}$, we next notice that the fact that $T \leq \frac{9p^2\pi^4}{4}\left(\sqrt{2}L\left(\sqrt{2}\|w_0\|_\infty + M\right)\sum_{n=1}^{\infty} n^{-4/3}\left(\max\left(2,\left(\sqrt{2}\|w_0\|_\infty + M\right)M^{-1/3}\right)\right)\right)^{-3}$ guarantees that $P$ maps $S$ into $S$. Indeed, for every $w\in S$, the fact that $F(0)=0$ and (3.1) imply the following inequalities:

$$\left| 2\int_0^t \exp(-pn^2\pi^2(t-\tau))\left(\int_0^1 (\bar{F}(w[\tau]))(s)\sin(n\pi s)ds\right)d\tau \right|$$

$$\leq \sqrt{2}\left(\sqrt{2}\|w_0\|_\infty + M\right)L\left(\sqrt{2}\|w_0\|_\infty + M\right)\int_0^t \exp(-pn^2\pi^2(t-\tau))d\tau$$

$$\leq \sqrt{2}\left(\sqrt{2}\|w_0\|_\infty + M\right)L\left(\sqrt{2}\|w_0\|_\infty + M\right)\exp(-pn^2\pi^2 t)\left(\int_0^t \exp(3pn^2\pi^2\tau/2)d\tau\right)^{2/3} t^{1/3},$$

$$\leq \sqrt{2}\left(\sqrt{2}\|w_0\|_\infty + M\right)L\left(\sqrt{2}\|w_0\|_\infty + M\right)\exp(-pn^2\pi^2 t)\left(2\frac{\exp(3pn^2\pi^2 t/2)}{3pn^2\pi^2}\right)^{2/3} t^{1/3}$$

$$\leq \sqrt{2}\left(\sqrt{2}\|w_0\|_\infty + M\right)L\left(\sqrt{2}\|w_0\|_\infty + M\right)\left(\frac{2}{3p\pi^2}\right)^{2/3} T^{1/3} n^{-4/3}$$

for all $t \in [0,T]$ and $n = 1, 2, \ldots$ (3.12)

Notice that we have used above Holder's inequality as well as the facts that $t \in [0,T]$ and $2\int_0^1 \sin^2(n\pi x)dx = 1$ for $n = 1, 2, \ldots$. Inequality (3.12) in conjunction with the fact that $T \leq \frac{9p^2\pi^4}{4}\left(\sqrt{2}L\left(\sqrt{2}\|w_0\|_\infty + M\right)\sum_{n=1}^\infty n^{-4/3}\left(\max\left(2,\left(\sqrt{2}\|w_0\|_\infty + M\right)M^{-1/3}\right)\right)\right)^{-3}$, definition (3.10) and estimate (3.9) guarantees that $Pw \in S$. Finally using a similar line of operations as in (3.12) we get for every $w, u \in S$:

$$\left| 2\int_0^t \exp(-pn^2\pi^2(t-\tau))\left(\int_0^1 (\bar{F}(w[\tau]) - \bar{F}(u[\tau]))(s)\sin(n\pi s)ds\right)d\tau \right|,$$

$$\leq \sqrt{2}L\left(\sqrt{2}\|w_0\|_\infty + M\right)\left(\frac{2}{3p\pi^2}\right)^{2/3} T^{1/3} n^{-4/3} \max_{0\leq\tau\leq T}\left(\|w[\tau] - u[\tau]\|_\infty\right)$$

for all $t \in [0,T]$ and $n = 1, 2, \ldots$ (3.13)

Inequality (3.13) in conjunction with the fact that $T \leq \frac{9p^2\pi^4}{4}\left(\sqrt{2}L\left(\sqrt{2}\|w_0\|_\infty + M\right)\sum_{n=1}^\infty n^{-4/3}\left(\max\left(2,\left(\sqrt{2}\|w_0\|_\infty + M\right)M^{-1/3}\right)\right)\right)^{-3}$, definition (3.10) guarantees that

$$\max_{0\leq t\leq T, 0\leq x\leq 1}\left(|(Pw)(t,x) - (Pu)(t,x)|\right) \leq \frac{1}{2}\max_{0\leq t\leq T, 0\leq x\leq 1}\left(|w(t,x) - u(t,x)|\right) \quad (3.14)$$

Therefore, the mapping $P: S \to S$ is a contraction. The conclusion of the lemma follows from Banach's fixed point theorem. The proof is complete. ◁

The following lemma guarantees important regularity properties for the solutions of (3.7).

**Lemma 3.3:** *Let $p > 0$ be a constant and let $\bar{F}: C^0([0,1]) \to C^0([0,1])$ be a continuous mapping with $\bar{F}(0) = 0$, for which there exists a non-decreasing function $L: \Re_+ \to (0, +\infty)$ such that (3.1) holds. Let $w_0 \in H_0^2(0,1)$ be given and let $w \in C^0([0,T] \times [0,1])$ be a solution of (3.7) for certain $T > 0$. Then $w \in C^1((0,T]; L^2(0,1))$ with $w[t] \in H_0^2(0,1)$ for all $t \in [0,T]$. Moreover, the mapping $t \to \|w_x[t]\|_2^2$ is $C^1$ on $[0,T]$ and the following equations hold:*

$$w_t[t] = pw_{xx}[t] + \bar{F}(w[t]), \text{ for } t \in (0,T] \quad (3.15)$$

$$\frac{d}{dt}\|w_x[t]\|_2^2 = -2p\|w_{xx}[t]\|_2^2 - 2\langle \bar{F}(w[t]), w_{xx}[t]\rangle, \text{ for } t \in [0,T] \tag{3.16}$$

**Proof:** Let $w \in C^0([0,T]\times[0,1])$ be a solution of (3.7). Notice that since

$$\left|\exp(-pn^2\pi^2 t) - \exp(-pn^2\pi^2 t_0)\right| \leq (pn^2\pi^2)^{1/r}\left(\frac{r-1}{r}\right)^{(r-1)/r} \exp(-pn^2\pi^2 t_0)(t-t_0)^{1/r}$$

for all $n \geq 1$, $t > t_0 \geq 0$, $r > 1$,

and since $\int_0^1 w_0(s)\sin(n\pi s)ds = -\frac{1}{n^2\pi^2}\int_0^1 w_0''(s)\sin(n\pi s)ds$ for $n=1,2,...$, we obtain from (3.7) for all $r > 2$, $t, t_0 \in [0,T]$ with $t > t_0$:

$$\begin{aligned}\|w[t]-w[t_0]\|_\infty &\leq 2\sum_{n=1}^\infty \left[\frac{1}{n^2\pi^2}\left|\exp(-pn^2\pi^2 t) - \exp(-pn^2\pi^2 t_0)\right|\int_0^1 |w_0''(s)\sin(n\pi s)|ds\right]\\ &+ 2\sum_{n=1}^\infty\left[\int_{t_0}^t \exp(-pn^2\pi^2(t-\tau))\left(\int_0^1|(\bar{F}(w[\tau]))(s)\sin(n\pi s)|ds\right)d\tau\right]\\ &+ 2\sum_{n=1}^\infty\left[\left|\exp(-pn^2\pi^2 t) - \exp(-pn^2\pi^2 t_0)\right|\int_0^{t_0}\exp(pn^2\pi^2 \tau)\left(\int_0^1 |(\bar{F}(w[\tau]))(s)\sin(n\pi s)|ds\right)d\tau\right]\\ &\leq \sqrt{2}\left(\frac{r-1}{pr\pi^2}\right)^{(r-1)/r}\left(p\|w_0''\|_2 + 2\max_{0\leq \tau\leq t}\left(\|\bar{F}(w[\tau])\|_2\right)\right)(t-t_0)^{1/r}\sum_{n=1}^\infty n^{-2(r-1)/r}\end{aligned} \tag{3.17}$$

Notice that for the above derivation we have used the Cauchy-Schwarz inequality and the facts that $|\sin(n\pi x)| \leq 1$ and $2\int_0^1 \sin^2(n\pi x)dx = 1$ for $n=1,2,...$. It follows from (3.17) that the following inequality for all $r > 2$, $t,t_0 \in [0,T]$:

$$\|w[t]-w[t_0]\|_\infty \leq K(T,r,p,w)|t-t_0|^{1/r} \tag{3.18}$$

where $K(T,r,p,w) := \sqrt{2}\left(\frac{r-1}{pr\pi^2}\right)^{(r-1)/r}\left(p\|w_0''\|_2 + 2\max_{0\leq t\leq T}\left(\|\bar{F}(w[t])\|_2\right)\right)\sum_{n=1}^\infty n^{-2(r-1)/r}$. Setting $\phi_n(x) = \sqrt{2}\sin(n\pi x)$ for $x \in [0,1]$, we get from (3.7) for $n=1,2,...$ and $t \in [0,T]$:

$$\begin{aligned}n^2\pi^2\langle w[t], \phi_n\rangle &= p^{-1}\left(1-\exp(-pn^2\pi^2 t)\right)\langle\bar{F}(w[t]),\phi_n\rangle - \exp(-pn^2\pi^2 t)\langle w_0'',\phi_n\rangle\\ &+ n^2\pi^2\int_0^t \exp(-pn^2\pi^2(t-\tau))\langle\bar{F}(w[\tau])-\bar{F}(w[t]),\phi_n\rangle d\tau\end{aligned} \tag{3.19}$$

Using (3.1), (3.18) and (3.19) we get for $n=1,2,...$, $r>2$ and $t \in [0,T]$:

$$\begin{aligned}n^2\pi^2|\langle w[t],\phi_n\rangle| &\leq p^{-1}|\langle\bar{F}(w[t]),\phi_n\rangle| + |\langle w_0'',\phi_n\rangle| + n^2\pi^2\int_0^t \exp(-pn^2\pi^2(t-\tau))|\langle\bar{F}(w[\tau])-\bar{F}(w[t]),\phi_n\rangle|d\tau\\ &\leq p^{-1}|\langle\bar{F}(w[t]),\phi_n\rangle| + |\langle w_0'',\phi_n\rangle| + n^2\pi^2 L\left(\max_{0\leq\tau\leq t}(\|w[\tau]\|_\infty)\right)\int_0^t \exp(-pn^2\pi^2(t-\tau))\|w[\tau]-w[t]\|_\infty d\tau\\ &\leq p^{-1}|\langle\bar{F}(w[t]),\phi_n\rangle| + |\langle w_0'',\phi_n\rangle| + n^2\pi^2 L\left(\max_{0\leq\tau\leq t}(\|w[\tau]\|_\infty)\right)K(T,r,p,w)\int_0^t (t-\tau)^{1/r}\exp(-pn^2\pi^2(t-\tau))d\tau\end{aligned} \tag{3.20}$$

We also have for $n=1,2,...$, $r \in (2,4)$, $a \in \left(\frac{r}{2}, 2\right)$ and $t \in [0,T]$:

$$\int_0^t (t-\tau)^{1/r} \exp\left(-pn^2\pi^2(t-\tau)\right) d\tau$$

$$\leq \int_{\max(0,t-n^{-a})}^{t} (t-\tau)^{1/r} \exp\left(-pn^2\pi^2(t-\tau)\right) d\tau + \int_0^{\max(0,t-n^{-a})} (t-\tau)^{1/r} \exp\left(-pn^2\pi^2(t-\tau)\right) d\tau$$

$$\leq n^{-a/r} \frac{1}{pn^2\pi^2} + T^{1/r} \frac{\exp\left(-pn^{2-a}\pi^2\right)}{pn^2\pi^2}$$

Using the fact that $n \exp\left(-pn^{2-a}\pi^2\right) \leq \left(\frac{1}{pe(2-a)\pi^2}\right)^{1/(2-a)}$, for $n=1,2,\ldots$, $a \in \left(\frac{r}{2}, 2\right)$, $r \in (2,4)$ in conjunction with the above inequality, we get for $n=1,2,\ldots$, $r \in (2,4)$, $a \in \left(\frac{r}{2}, 2\right)$ and $t \in [0,T]$:

$$\int_0^t (t-\tau)^{1/r} \exp\left(-pn^2\pi^2(t-\tau)\right) d\tau \leq \frac{n^{-a/r}}{pn^2\pi^2}\left(1 + T^{1/r}\left(\frac{1}{p\pi^2 e(2-a)}\right)^{1/(2-a)}\right) \quad (3.21)$$

Combining (3.21) with (3.20) we get for $n=1,2,\ldots$, $r \in (2,4)$, $a \in \left(\frac{r}{2}, 2\right)$ and $t \in [0,T]$:

$$n^2\pi^2 |\langle w[t], \phi_n \rangle| \leq p^{-1}|\langle \overline{F}(w[t]), \phi_n \rangle| + |\langle w_0'', \phi_n \rangle|$$
$$+ n^{-a/r} p^{-1} L\left(\max_{0 \leq \tau \leq T}\left(\|w[\tau]\|_\infty\right)\right) K(T,r,p,w)\left(1 + T^{1/r}\left(\frac{1}{p\pi^2 e(2-a)}\right)^{1/(2-a)}\right) \quad (3.22)$$

Inequality (3.22) in conjunction with the fact that $\{\phi_n : n=1,2,\ldots\}$ is an orthonormal basis of $L^2(0,1)$ and Parseval's identity, shows that $w[t] \in H_0^2(0,1)$ for all $t \in [0,T]$ with

$$w_{xx}[t] = \sum_{n=1}^{\infty} \phi_n \exp\left(-pn^2\pi^2 t\right)\langle w_0'', \phi_n \rangle - p^{-1} \sum_{n=1}^{\infty} \phi_n \left(1 - \exp\left(-pn^2\pi^2 t\right)\right)\langle \overline{F}(w[t]), \phi_n \rangle$$
$$- \sum_{n=1}^{\infty} n^2\pi^2 \phi_n \int_0^t \exp\left(-pn^2\pi^2(t-\tau)\right)\langle \overline{F}(w[\tau]) - \overline{F}(w[t]), \phi_n \rangle d\tau \quad (3.23)$$

$$w_x[t] = p^{-1} \sum_{n=1}^{\infty} n^{-1}\pi^{-1} \psi_n \left(1 - \exp\left(-pn^2\pi^2 t\right)\right)\langle \overline{F}(w[t]), \phi_n \rangle$$
$$- \sum_{n=1}^{\infty} n^{-1}\pi^{-1} \psi_n \exp\left(-pn^2\pi^2 t\right)\langle w_0'', \phi_n \rangle \quad (3.24)$$
$$+ \sum_{n=1}^{\infty} n\pi \psi_n \int_0^t \exp\left(-pn^2\pi^2(t-\tau)\right)\langle \overline{F}(w[\tau]) - \overline{F}(w[t]), \phi_n \rangle d\tau$$

for $t \in [0,T]$, where $\psi_n(x) = \sqrt{2}\cos(n\pi x)$ for $x \in [0,1]$. Differentiating formally with respect to $t \in [0,T]$ the right hand side of equation (3.7) and using the fact that $\int_0^1 w_0(s)\sin(n\pi s)ds = -\frac{1}{n^2\pi^2}\int_0^1 w_0''(s)\sin(n\pi s)ds$ for $n=1,2,\ldots$, we obtain the Fourier series:

$$p\sum_{n=1}^{\infty} \phi_n \exp\left(-pn^2\pi^2 t\right)\langle w_0'', \phi_n \rangle + \sum_{n=1}^{\infty} \phi_n \exp\left(-pn^2\pi^2 t\right)\langle \overline{F}(w[t]), \phi_n \rangle$$
$$- p\sum_{n=1}^{\infty} \phi_n n^2\pi^2 \int_0^t \exp\left(-pn^2\pi^2(t-\tau)\right)\langle \overline{F}(w[\tau]) - \overline{F}(w[t]), \phi_n \rangle d\tau$$

Similarly as above and using the fact that $pn^2\pi^2 t \exp\left(-pn^2\pi^2 t\right) \leq e^{-1}$, we show that the absolute value of each Fourier coefficient of the above series is bounded by

$$n^{-2}\pi^{-2}t_0^{-1}e^{-1}p^{-1}\max_{0\leq\tau\leq T}\left(\left\|\bar{F}(w[\tau])\right\|_2\right)+p\left|\langle w_0'',\phi_n\rangle\right|$$

$$+n^{-a/r}L\left(\max_{0\leq\tau\leq t}\left(\left\|w[\tau]\right\|_\infty\right)\right)K(T,r,p,w)\left(1+T^{1/r}\left(\frac{1}{p\pi^2 e(2-a)}\right)^{1/(2-a)}\right)$$

for $n=1,2,...$, $r\in(2,4)$, $a\in\left(\frac{r}{2},2\right)$, $t_0\in(0,T)$ and $t\in[t_0,T]$. The above inequality shows that $w\in C^1\left((0,T];L^2(0,1)\right)$ with

$$w_t[t]=p\sum_{n=1}^\infty \phi_n \exp\left(-pn^2\pi^2 t\right)\langle w_0'',\phi_n\rangle + \sum_{n=1}^\infty \phi_n \exp\left(-pn^2\pi^2 t\right)\langle \bar{F}(w[t]),\phi_n\rangle$$
$$-p\sum_{n=1}^\infty \phi_n n^2\pi^2 \int_0^t \exp\left(-pn^2\pi^2(t-\tau)\right)\langle \bar{F}(w[\tau])-\bar{F}(w[t]),\phi_n\rangle d\tau \quad (3.25)$$

for all $t\in(0,T]$. Equation (3.15) is a direct consequence of equations (3.23) and (3.25).

Finally, we show the validity of equation (3.16). Writing (3.7), (3.23) and (3.24) in a different way (using the fact that $\int_0^1 w_0(s)\sin(n\pi s)ds = -\frac{1}{n^2\pi^2}\int_0^1 w_0''(s)\sin(n\pi s)ds$ for $n=1,2,...$), we get for all $t\in[0,T]$:

$$w_{xx}[t]=\sum_{n=1}^\infty \phi_n\left[\exp\left(-pn^2\pi^2 t\right)\langle w_0'',\phi_n\rangle - n^2\pi^2\int_0^t \exp\left(-pn^2\pi^2(t-\tau)\right)\langle \bar{F}(w[\tau]),\phi_n\rangle d\tau\right] \quad (3.26)$$

$$w_x[t]=\sum_{n=1}^\infty \psi_n\left[-n^{-1}\pi^{-1}\exp\left(-pn^2\pi^2 t\right)\langle w_0'',\phi_n\rangle + n\pi\int_0^t \exp\left(-pn^2\pi^2(t-\tau)\right)\langle \bar{F}(w[\tau]),\phi_n\rangle d\tau\right] \quad (3.27)$$

$$w[t]=\sum_{n=1}^\infty \phi_n\left[-n^{-2}\pi^{-2}\exp\left(-pn^2\pi^2 t\right)\langle w_0'',\phi_n\rangle + \int_0^t \exp\left(-pn^2\pi^2(t-\tau)\right)\langle \bar{F}(w[\tau]),\phi_n\rangle d\tau\right] \quad (3.28)$$

$$\left\|w_x[t]\right\|_2^2 = \sum_{n=1}^\infty \left[n^{-1}\pi^{-1}\exp\left(-pn^2\pi^2 t\right)\langle w_0'',\phi_n\rangle - n\pi\int_0^t \exp\left(-pn^2\pi^2(t-\tau)\right)\langle \bar{F}(w[\tau]),\phi_n\rangle d\tau\right]^2 \quad (3.29)$$

Notice that (3.29) is derived from (3.27) by using the fact that $\{\psi_n: n=0,1,2,...\}$ with $\psi_0(x)=1$ for $x\in[0,1]$ is an orthonormal basis of $L^2(0,1)$ in conjunction with Parseval's identity. Differentiating formally the right hand side of (3.29) with respect to $t\in[0,T]$, we get (using (3.26), (3.27), (3.28)) the following series for all $t\in[0,T]$:

$$-2p\sum_{n=1}^\infty n^4\pi^4 \langle w[t],\phi_n\rangle^2 + 2\sum_{n=1}^\infty n^2\pi^2 \langle \bar{F}(w[t]),\phi_n\rangle\langle w[t],\phi_n\rangle$$

which by virtue of (3.22) converges uniformly and absolutely for all $t\in[0,T]$. Using (3.26) and (3.28), we obtain (3.16) for all $t\in[0,T]$. The proof is complete. ◁

The following lemma clarifies what happens when the solution of (3.7) cannot be continued.

**Lemma 3.4:** *Let $p>0$ be a constant and let $\bar{F}:C^0([0,1])\to C^0([0,1])$ be a continuous mapping with $\bar{F}(0)=0$, for which there exists a non-decreasing function $L:\Re_+\to(0,+\infty)$ such that (3.1) holds.*

*Then for every $w_0 \in H_0^2(0,1)$ there exists $t_{\max}(w_0) \in (0,+\infty]$ and a mapping $w \in C^0([0,t_{\max}(w_0)) \times [0,1])$ such that (3.7) holds for all $(t,x) \in [0,t_{\max}(w_0)) \times [0,1]$. Moreover, if $t_{\max}(w_0) < +\infty$, then $\lim_{t \to t_{\max}^-} (\|w[t]\|_\infty) = +\infty$.*

**Proof:** Define for each $w_0 \in H_0^2(0,1)$

$$t_{\max}(w_0) := \sup\{T > 0 : \exists w \in C^0([0,T] \times [0,1]) \text{ that satisfies (33)}\} \tag{3.30}$$

and notice that the set $\{T > 0 : \exists w \in C^0([0,T] \times [0,1]) \text{ that satisfies (33)}\}$ is non-empty (by virtue of Lemma 3.2). Suppose that $t_{\max}(w_0) < +\infty$. For every $T \in (0, t_{\max}(w_0))$ there exists $w \in C^0([0,T] \times [0,1])$ that satisfies (3.7). Consider a solution $u \in C^0([0,T'] \times [0,1])$ of

$$\begin{aligned} u(t,x) &= 2\sum_{n=1}^\infty \left[ \exp(-pn^2\pi^2 t) \int_0^1 u_0(s) \sin(n\pi s) ds \right] \sin(n\pi x) \\ &+ 2\sum_{n=1}^\infty \left[ \int_0^t \exp(-pn^2\pi^2 (t-\tau)) \left( \int_0^1 (\bar{F}(u[\tau]))(s) \sin(n\pi s) ds \right) d\tau \right] \sin(n\pi x) \end{aligned} \tag{3.31}$$

where $u_0 = w[T]$. Notice that Lemma 3.3 guarantees that $u_0 \in H_0^2(0,1)$. By virtue of Lemma 3.2, such a solution $u \in C^0([0,T'] \times [0,1])$ exists for all $T' > 0$ with $T' \leq \frac{9p^2\pi^4}{8\sqrt{2}} \left( \left( \max(2, \sqrt{2}\|w[T]\|_\infty + 1) \right) L\left( \sqrt{2}\|w[T]\|_\infty + 1 \right) \sum_{n=1}^\infty n^{-4/3} \right)^{-3}$. Moreover, consider the function:

$$\tilde{w}[t] := \begin{cases} w[t] & \text{if } t \in [0,T] \\ u[t-T] & \text{if } t \in (T, T+T'] \end{cases} \tag{3.32}$$

Notice that $\tilde{w} \in C^0([0, T+T'] \times [0,1])$ and that $\tilde{w}$ is a solution of (3.7). Thus, definition (3.30) implies that

$$\frac{9p^2\pi^4}{8\sqrt{2}} \left( \left( \max(2, \sqrt{2}\|w[T]\|_\infty + 1) \right) L\left( \sqrt{2}\|w[T]\|_\infty + 1 \right) \sum_{n=1}^\infty n^{-4/3} \right)^{-3} + T \leq t_{\max}(w_0)$$

which shows that $\lim_{T \to t_{\max}^-} (\|w[T]\|_\infty) = +\infty$. The proof is complete. ◁

The following theorem deals with a similar case to that of Theorem 3.1 but the boundary conditions are homogeneous.

**Theorem 3.5:** *Let $p > 0$ be a constant and let $\bar{F} : C^0([0,1]) \to C^0([0,1])$ be a continuous mapping with $\bar{F}(0) = 0$, for which there exists a non-decreasing function $L : \Re_+ \to (0,+\infty)$ such that (3.1) holds. Then for every $w_0 \in H_0^2(0,1)$ there exists $t_{\max}(w_0) \in (0,+\infty]$ and a unique mapping $w \in C^0([0, t_{\max}(w_0)) \times [0,1]) \cap C^1((0, t_{\max}(w_0)); L^2(0,1))$ with $w[0] = w_0$, $w[t] \in H_0^2(0,1)$ for all $t \in [0, t_{\max}(w_0))$, for which the mapping $t \to \|w_x[t]\|_2^2$ is $C^1$ on $[0, t_{\max}(w_0))$ and for which the following equations hold:*

$$w_t[t] = p w_{xx}[t] + \bar{F}(w[t]), \text{ for } t \in (0, t_{\max}(w_0)) \tag{3.33}$$

$$\frac{d}{dt} \|w_x[t]\|_2^2 = -2p \|w_{xx}[t]\|_2^2 - 2 \langle \bar{F}(w[t]), w_{xx}[t] \rangle, \text{ for } t \in [0, t_{\max}(w_0)) \tag{3.34}$$

*Moreover, if $t_{\max}(w_0) < +\infty$, then $\lim_{t \to t_{\max}^-} (\|w[t]\|_\infty) = +\infty$.*

**Proof:** The result is a direct consequence of Lemma 3.3 and Lemma 3.4. The only thing that needs to be shown is uniqueness.

Let $w_0 \in H_0^2(0,1)$ and $T > 0$ be given and consider two solutions $w, u \in C^0([0,T] \times [0,1]) \cap C^1((0,T]; L^2(0,1))$ with $w[0] = u[0] = w_0$, $w[t], u[t] \in H_0^2(0,1)$ for all $t \in [0,T]$, for which the equations $w_t[t] = pw_{xx}[t] + \bar{F}(w[t])$ and $u_t[t] = pu_{xx}[t] + \bar{F}(u[t])$ hold for $t \in [0,T]$. Notice that necessarily $w, u$ satisfy

$$w(t,x) = 2\sum_{n=1}^{\infty} \left[ \exp\left(-pn^2\pi^2(t-t_0)\right) \int_0^1 (w[t_0])(s) \sin(n\pi s) ds \right] \sin(n\pi x)$$
$$+ 2\sum_{n=1}^{\infty} \left[ \int_0^{t-t_0} \exp\left(-pn^2\pi^2(t-t_0-\tau)\right) \left( \int_0^1 (\bar{F}(w[t_0+\tau]))(s) \sin(n\pi s) ds \right) d\tau \right] \sin(n\pi x) \quad (3.35)$$

$$u(t,x) = 2\sum_{n=1}^{\infty} \left[ \exp\left(-pn^2\pi^2(t-t_0)\right) \int_0^1 (u[t_0])(s) \sin(n\pi s) ds \right] \sin(n\pi x)$$
$$+ 2\sum_{n=1}^{\infty} \left[ \int_0^{t-t_0} \exp\left(-pn^2\pi^2(t-t_0-\tau)\right) \left( \int_0^1 (\bar{F}(u[t_0+\tau]))(s) \sin(n\pi s) ds \right) d\tau \right] \sin(n\pi x) \quad (3.36)$$

for all $t_0, x \in [0,T] \times [0,1]$ and $t \in [t_0, T]$. Notice that (3.35), (3.36) give (3.7) and (3.31) for $t_0 = 0$. Define

$$M := \max\left( \max_{0 \le t \le T} \left( \|w[t]\|_{\infty} \right), \max_{0 \le t \le T} \left( \|u[t]\|_{\infty} \right) \right) \quad (3.37)$$

$$r = \frac{9p^2\pi^4}{4} \left( \sqrt{2}L\left((1+\sqrt{2})M\right) \sum_{n=1}^{\infty} n^{-4/3} \left( \max\left(2, \left((1+\sqrt{2})M\right)M^{-1/3}\right) \right) \right)^{-3} \quad (3.38)$$

By virtue of Lemma 3.2 and the fact that $w[0] = u[0] = w_0$, we have $w[t] = u[t]$ on $[0, \min(r,T)]$. If $r \ge T$ then we conclude that $w = u$ on $[0,T]$. On the other hand, if $r < T$ then we apply Lemma 3.2 to equations (3.35), (3.36) with $t_0 = r$. Using the fact that $w[r] = u[r]$, we get $w[t] = u[t]$ on $[0, \min(2r,T)]$. Similarly, we may continue to show that $w[t] = u[t]$ on $[0, \min(nr,T)]$ for all $n \ge 1$. Thus, we conclude that $w = u$ on $[0,T]$.

Since $T > 0$ is arbitrary, we have shown uniqueness of the mapping $w \in C^0([0, t_{\max}(w_0)) \times [0,1]) \cap C^1((0, t_{\max}(w_0)); L^2(0,1))$ with $w[0] = w_0$, $w[t] \in H_0^2(0,1)$ for all $t \in [0, t_{\max}(w_0))$, which satisfies (3.33). The proof is complete. ◁

We are now ready to provide the proof of Theorem 3.1.

**Proof of Theorem 3.1:** Consider the Volterra operator $K : L^2(0,1) \to L^2(0,1)$ defined by (3.5) with inverse

$$(K^{-1}w)(x) = w(x) + \int_0^x \bar{k}(\tau) \exp\left( \int_\tau^x \bar{k}(s) ds \right) w(\tau) d\tau \text{ for } w \in L^2(0,1), \ x \in [0,1] \quad (3.39)$$

Consider also the linear continuous operator $G : L^2(0,1) \to L^2(0,1)$ defined by (3.6). Using (3.5), (3.6), we can guarantee that the transformation $w[t] = Ku[t]$ transforms the solutions of (3.2), (3.3) with initial condition $u[0] = u_0$ to the solutions of

$$w_t[t] = pw_{xx}[t] + pGK^{-1}w[t] + K\bar{F}(K^{-1}w[t]) \quad (3.40)$$

$$w(t,0) = w(t,1) = 0 \quad (3.41)$$

with initial condition $w[0] = w_0$ given by

$$w_0 = Ku_0 \quad (3.42)$$

Moreover, using (3.5), (3.6), we can guarantee that the transformation $u[t] = K^{-1}w[t]$ transforms the solutions of (3.40), (3.41) with initial condition $w[0] = w_0$ to the solutions of (3.2), (3.3) with $u_0 = K^{-1}u_0$.

Notice that by virtue of (3.1) and the fact that both $K, K^{-1}$ are continuous linear operators on $C^0([0,1])$, the (nonlinear) operator $\tilde{F}: C^0([0,1]) \to C^0([0,1])$ defined by $\tilde{F}(w) = pGK^{-1}w + K\bar{F}(K^{-1}w)$ for all $w \in C^0([0,1])$ satisfies (3.1) with $\bar{F}$ replaced by $\tilde{F}$ (and a different function $L$). Theorem 3.5 implies that the solution of (3.40), (3.41) with initial condition given by (3.42) for $u_0 \in H^2(0,1)$ with $u_0(0) = 0$, $u_0(1) = \langle \bar{k}, u_0 \rangle$, exists locally and is unique. Therefore, the solution of (3.2), (3.3) with initial condition $u[0] = u_0$ exists locally and is unique. Inequality (3.4) follows directly from (3.34), (3.5) and the fact that $w[t] = Ku[t]$.
The proof is complete. ◁

**Remark 3.6:** The proof of Theorem 3.1, which uses Theorem 3.5, shows one important technical issue that would arise in the case of a nonlinear boundary feedback law: since the solution given by Theorem 3.5 is of class $C^1((0, t_{\max}(w_0)); L^2(0,1))$, the time derivative of the control action would not be differentiable if the feedback law is nonlinear. That would not allow the homogenization of the boundary conditions (as performed in the proof of Theorem 3.1), which is necessary for the development of existence/uniqueness results.

## 4. Proof of Main Result

For the proof of Theorem 2.1, several auxiliary results are needed. The auxiliary results are of independent interest. We start with the following proposition which extends the well-known Wirtinger's inequality.

**Proposition 4.1 (Extension of Wirtinger's inequality):** *For every $u \in H_0^1(0,1)$ and for every non-zero $k \in L^2(0,1)$ the following inequality holds for every $\varepsilon \geq 0$:*

$$\left(4 - \frac{3(1+\varepsilon)\left(\|k\|_2^2 - k_1^2\right)}{(1+\varepsilon)\|k\|_2^2 - k_1^2}\right)\pi^2 \|u\|_2^2 \leq \|u'\|_2^2 + \frac{3\pi^2 \varepsilon(1+\varepsilon)}{(1+\varepsilon)\|k\|_2^2 - k_1^2}\langle k, u \rangle^2 \quad (4.1)$$

*where* $k_1 = \sqrt{2}\int_0^1 k(x)\sin(\pi x)dx$.

**Remark 4.2:** For $\varepsilon = 0$ inequality (4.1) gives the well-known Wirtinger's inequality: $\pi^2 \|u\|_2^2 \leq \|u'\|_2^2$.

**Proof:** First notice that when $k_1 = 0$ inequality (4.1) becomes

$$\pi^2 \|u\|_2^2 \leq \|u'\|_2^2 + \frac{3\pi^2 \varepsilon}{\|k\|_2^2}\langle k, u \rangle^2$$

which holds for all $\varepsilon \geq 0$ (by virtue of Wirtinger's inequality). Similarly, for $\varepsilon = 0$ inequality (4.1) gives Wirtinger's inequality.

Next we assume that $k_1 \neq 0$ and we show inequality (4.1) for every $\varepsilon > 0$. Since $u \in H^1(0,1)$ with $u(0) = u(1) = 0$, it follows that

$$u = \sum_{n=1}^{\infty} \langle u, \phi_n \rangle \phi_n \quad (4.2)$$

$$u' = \sum_{n=1}^{\infty} n\pi \langle u, \phi_n \rangle \psi_n \tag{4.3}$$

where $\phi_n(x) = \sqrt{2}\sin(n\pi x)$, $\psi_n(x) = \sqrt{2}\cos(n\pi x)$ and the above equalities are in the sense of $L^2(0,1)$. Since $k \in L^2(0,1)$, we get from (4.2) and definition $k_1 = \sqrt{2}\int_0^1 k(x)\sin(\pi x)dx$:

$$\langle k, u \rangle = k_1 \langle u, \phi_1 \rangle + \sum_{n=2}^{\infty} \langle k, \phi_n \rangle \langle u, \phi_n \rangle \tag{4.4}$$

Moreover, since $\{\phi_n : n = 1, 2, ...\}$ and $\{\psi_n : n = 0, 1, 2, ...\}$ with $\psi_0(x) = 1$ for $x \in [0,1]$ are orthonormal bases of $L^2(0,1)$, we obtain from (4.2), (4.3) and Parseval's identity:

$$\|u\|_2^2 = \sum_{n=1}^{\infty} \langle u, \phi_n \rangle^2, \quad \|k\|_2^2 = \sum_{n=1}^{\infty} \langle k, \phi_n \rangle^2 \text{ and } \|u'\|_2^2 = \sum_{n=1}^{\infty} n^2 \pi^2 \langle u, \phi_n \rangle^2 \tag{4.5}$$

Using (4.4) and the Cauchy-Schwarz inequality, we get:

$$\left| \langle u, \phi_1 \rangle - k_1^{-1} \langle k, u \rangle \right| \leq |k_1|^{-1} \left( \sum_{n=2}^{\infty} \langle k, \phi_n \rangle^2 \right)^{1/2} \left( \sum_{n=2}^{\infty} \langle u, \phi_n \rangle^2 \right)^{1/2} \tag{4.6}$$

The triangle inequality (which implies that $|\langle u, \phi_1 \rangle| \leq |\langle u, \phi_1 \rangle - k_1^{-1}\langle k, u \rangle| + |k_1^{-1}\langle k, u \rangle|$) gives the following inequality for all $\varepsilon > 0$

$$\left| \langle u, \phi_1 \rangle \right|^2 \leq (1 + \varepsilon^{-1}) \left| \langle u, \phi_1 \rangle - k_1^{-1} \langle k, u \rangle \right|^2 + (1 + \varepsilon) \left| k_1^{-1} \langle k, u \rangle \right|^2$$

which in conjunction with (4.5) and (4.6) (that imply $\left| \langle u, \phi_1 \rangle - k_1^{-1} \langle k, u \rangle \right| \leq |k_1|^{-1} \left( \|k\|_2^2 - k_1^2 \right)^{1/2} \left( \|u\|_2^2 - \langle u, \phi_1 \rangle^2 \right)^{1/2}$) gives for all $\varepsilon > 0$:

$$\left| \langle u, \phi_1 \rangle \right|^2 \leq \frac{(1+\varepsilon)\left( \|k\|_2^2 - k_1^2 \right)}{(1+\varepsilon)\|k\|_2^2 - k_1^2} \|u\|_2^2 + \frac{\varepsilon(1+\varepsilon)}{(1+\varepsilon)\|k\|_2^2 - k_1^2} \left| \langle k, u \rangle \right|^2 \tag{4.7}$$

Moreover, we obtain from (4.5):

$$\|u'\|_2^2 = \pi^2 \left| \langle u, \phi_1 \rangle \right|^2 + \sum_{n=2}^{\infty} n^2\pi^2 \langle u, \phi_n \rangle^2 \geq \pi^2 \left| \langle u, \phi_1 \rangle \right|^2 + 4\pi^2 \sum_{n=2}^{\infty} \langle u, \phi_n \rangle^2 = 4\pi^2 \|u\|_2^2 - 3\pi^2 \left| \langle u, \phi_1 \rangle \right|^2$$

which implies the inequality

$$4\pi^2 \|u\|_2^2 \leq \|u'\|_2^2 + 3\pi^2 \left| \langle u, \phi_1 \rangle \right|^2 \tag{4.8}$$

Combining (4.7) and (4.8) and rearranging, we obtain (4.1) for all $\varepsilon > 0$.
The proof is complete. ◁

**Proposition 4.3:** *Suppose that $u \in H^1(0,1)$ satisfies $u(0) = 0$. Then for every $k \in H^1(0,1)$ with $k(0) = 0$ and $k(1) = 1$ the following inequality holds for every $\varepsilon \geq 0$:*

$$\frac{(1+\varepsilon)\|k\|_2^2 + (3\varepsilon - 1)k_1^2}{(1+\varepsilon)\|k\|_2^2 - k_1^2} \pi^2 \|u\|_2^2 \leq \|u'\|_2^2 + \frac{3\pi^2 \varepsilon(1+\varepsilon)}{(1+\varepsilon)\|k\|_2^2 - k_1^2} \langle k, u \rangle^2$$
$$+ u^2(1)\left( \|k'\|_2^2 + (3\varepsilon - 1)\pi^2 \|k\|_2^2 \right) - 2u(1)\left( (3\varepsilon - 1)\pi^2 \langle k, u \rangle + \langle k', u' \rangle \right) \tag{4.9}$$

*where $k_1 = \sqrt{2}\int_0^1 k(x)\sin(\pi x)dx$.*

**Proof:** Define
$$w := u - u(1)k \tag{4.10}$$

Since $u, k \in H^1(0,1)$ with $u(0) = k(0) = 0$ and $k(1) = 1$, it follows from definition (4.10) that $w \in H^1(0,1)$ satisfies $w(0) = w(1) = 0$. Therefore, by virtue of Proposition 4.1, it follows that the following inequality holds for every $\varepsilon \geq 0$:

$$\left(4 - \frac{3(1+\varepsilon)\left(\|k\|_2^2 - k_1^2\right)}{(1+\varepsilon)\|k\|_2^2 - k_1^2}\right)\pi^2 \|w\|_2^2 \leq \|w'\|_2^2 + \frac{3\pi^2 \varepsilon(1+\varepsilon)}{(1+\varepsilon)\|k\|_2^2 - k_1^2}\langle k, w\rangle^2 \tag{4.11}$$

where $k_1 = \sqrt{2}\int_0^1 k(x)\sin(\pi x)dx$. Notice that definition (4.10) implies the following equalities:

$$\begin{aligned}
\|w\|_2^2 &= \|u\|_2^2 + u^2(1)\|k\|_2^2 - 2u(1)\langle k, u\rangle \\
\|w'\|_2^2 &= \|u'\|_2^2 + u^2(1)\|k'\|_2^2 - 2u(1)\langle k', u'\rangle \\
\langle k, w\rangle &= \langle k, u\rangle - u(1)\|k\|_2^2
\end{aligned} \tag{4.12}$$

Combining (4.11) and (4.12), we obtain (4.9). The proof is complete. ◁

The following proposition uses the CLF defined by (2.12) as well as Proposition 4.3 in order to obtain a bound for the $L^2$ spatial norm of the state. It is one of the main auxiliary results.

**Proposition 4.4:** Let $t_{\max} \in (0, +\infty]$ and $u \in C^0\left([0, t_{\max}) \times [0,1]\right) \cap C^1\left((0, t_{\max}); L^2(0,1)\right)$ with $u[t] \in H^2(0,1)$ for all $t \in [0, t_{\max})$ be a function that satisfies (2.1) for $t \in (0, t_{\max})$ and (2.2), (2.6) for $t \in [0, t_{\max})$, where $k \in C^2([0,1])$ is a function with $k(0) = 0$, $k(1) = 1$ that satisfies $k''(x) = \mu k(x)$ for $x \in [0,1]$ and for some constant $\mu \in \Re$, $p > 0$, $r \in \Re$ are constants with $1 + \|k\|_2^2 r > 0$, and the mapping $F : C^0([0,1]) \to C^0([0,1])$ is defined by (2.3) for a locally Lipschitz function $f \in C^0([0,1] \times \Re)$ that satisfies (2.4), (2.5) for certain real constants $q \in \Re$, $\gamma, \delta, B \geq 0$ and $b \geq 2$ with $B \geq \delta |r| \|k\|_b \|k\|_{\frac{b}{b-1}}$. Suppose that there exists a constant $\varepsilon \geq 0$ that satisfies (2.7). Then there exist constants $G, \sigma > 0$ (independent of the particular solution $u$ and the time $t_{\max} \in (0, +\infty]$) such that (2.8) holds for all $t \in [0, t_{\max})$.

**Proof:** Consider the Lyapunov functional defined by (2.12), which (by virtue of the Cauchy-Schwarz inequality) satisfies for all $u \in L^2(0,1)$

$$\frac{1}{2}\left(1 + \min(0, r)\|k\|_2^2\right)\|u\|_2^2 \leq V(u) \leq \frac{1}{2}\left(1 + \max(0, r)\|k\|_2^2\right)\|u\|_2^2 \tag{4.13}$$

Therefore, by virtue of (4.13) and the inequality $1 + \|k\|_2^2 r > 0$, there exist constants $c_2 \geq c_1 > 0$ such that

$$c_1 \|u\|_2^2 \leq V(u) \leq c_2 \|u\|_2^2, \text{ for all } u \in L^2(0,1) \tag{4.14}$$

Since $u \in C^0\left([0, t_{\max}) \times [0,1]\right) \cap C^1\left((0, t_{\max}); L^2(0,1)\right)$, it follows from definition (2.12) that the mapping $(0, t_{\max}) \ni t \to V(u[t]) \in \Re$ is continuously differentiable with $V(u[t])$ being continuous at $t = 0$ and satisfies

$$\frac{d}{dt}V(u[t]) = \langle u[t], u_t[t]\rangle + r\langle k, u[t]\rangle\langle k, u_t[t]\rangle, \text{ for all } t \in (0, t_{\max}) \tag{4.15}$$

Using (2.1), (2.3), (2.4), the fact that $B \geq \delta |r| \|k\|_b \|k\|_{\frac{b}{b-1}}$, integration by parts and (2.6), we get from (4.15) for all $t \in (0, t_{max})$:

$$\frac{d}{dt} V(u[t]) \leq pu(t,1) u_x(t,1) - p \|u_x[t]\|_2^2 + (q+\gamma) \|u[t]\|_2^2 - \delta |r| \|k\|_b \|k\|_{\frac{b}{b-1}} \|u[t]\|_b^b$$
$$+ rpk(1)\langle k, u[t]\rangle u_x(t,1) - rp\langle k, u[t]\rangle \langle k', u_x[t]\rangle + rq\langle k, u[t]\rangle^2 \quad (4.16)$$
$$+ |r| |\langle k, u[t]\rangle| |\langle k, F(u[t]) - qu[t]\rangle|$$

Using (2.5), the Cauchy-Schwarz inequality and the Holder inequality, we get $|\langle k, F(u[t]) - qu[t]\rangle| \leq \gamma \|k\|_2 \|u[t]\|_2 + \delta \|k\|_b \|u[t]\|_b^{b-1}$. Combining the previous estimate with (4.16) and applying again integration by parts, we get for all $t \in (0, t_{max})$:

$$\frac{d}{dt} V(u[t]) \leq p\big(u(t,1) + rk(1)\langle k, u[t]\rangle\big) u_x(t,1) - p\|u_x[t]\|_2^2$$
$$+ (q+\gamma)\|u[t]\|_2^2 - \delta |r| \|k\|_b \|k\|_{\frac{b}{b-1}} \|u[t]\|_b^b$$
$$- rpk'(1) u(t,1) \langle k, u[t]\rangle + rp\langle k, u[t]\rangle \langle k'', u[t]\rangle + rq\langle k, u[t]\rangle^2 \quad (4.17)$$
$$+ \gamma |r| |\langle k, u[t]\rangle| \|k\|_2 \|u[t]\|_2 + \delta |r| |\langle k, u[t]\rangle| \|k\|_b \|u[t]\|_b^{b-1}$$

The Cauchy-Schwarz inequality implies that $|\langle k, u[t]\rangle| \|k\|_2 \|u[t]\|_2 \leq \|k\|_2^2 \|u[t]\|_2^2$. Moreover, using Holder's inequality we get $|\langle k, u[t]\rangle| \leq \|k\|_{\frac{b}{b-1}} \|u[t]\|_b$. Combining the previous inequalities with (4.17) and using (2.2), (2.6) and the fact that $k(1)=1$ and $k'' = \mu k$, we get for all $t \in (0, t_{max})$:

$$\frac{d}{dt} V(u[t]) \leq -p \|u_x[t]\|_2^2 + \left(q + \gamma\big(1 + |r| \|k\|_2^2\big)\right) \|u[t]\|_2^2$$
$$+ \big(r^2 pk'(1) + rp\mu + rq\big) \langle k, u[t]\rangle^2 \quad (4.18)$$

Exploiting Proposition 4.3 in conjunction with (4.18), we get the following inequality for every $\varepsilon \geq 0$ and $t \in (0, t_{max})$:

$$\frac{d}{dt} V(u[t]) \leq \left(q + \gamma\big(1 + |r|\|k\|_2^2\big) - p \frac{(1+\varepsilon)\|k\|_2^2 + (3\varepsilon - 1)k_1^2}{(1+\varepsilon)\|k\|_2^2 - k_1^2} \pi^2 \right) \|u[t]\|_2^2$$
$$+ p \frac{3\pi^2 \varepsilon(1+\varepsilon)}{(1+\varepsilon)\|k\|_2^2 - k_1^2} \langle k, u[t]\rangle^2 + \big(pr^2 k'(1) + pr\mu + rq\big) \langle k, u[t]\rangle^2 \quad (4.19)$$
$$+ pu^2(t,1)\big(\|k'\|_2^2 + (3\varepsilon - 1)\pi^2 \|k\|_2^2\big) - 2pu(t,1)\big((3\varepsilon - 1)\pi^2 \langle k, u[t]\rangle + \langle k', u_x[t]\rangle\big)$$

Using (2.6), the fact that $k'' = \mu k$, integration by parts and (4.19), we get for every $\varepsilon \geq 0$ and $t \in (0, t_{max})$:

$$\frac{d}{dt} V(u[t]) \leq \left(q + \gamma\big(1 + |r|\|k\|_2^2\big) - p \frac{(1+\varepsilon)\|k\|_2^2 + (3\varepsilon - 1)k_1^2}{(1+\varepsilon)\|k\|_2^2 - k_1^2} \pi^2 \right) \|u[t]\|_2^2$$
$$+ p\left(r^2\big(\|k'\|_2^2 + (3\varepsilon - 1)\pi^2\|k\|_2^2 - k'(1)\big) + r\left(2(3\varepsilon - 1)\pi^2 + \frac{q}{p} - \mu\right) + \frac{3\pi^2 \varepsilon(1+\varepsilon)}{(1+\varepsilon)\|k\|_2^2 - k_1^2}\right) \langle k, u[t]\rangle^2 \quad (4.20)$$

Let $\varepsilon \geq 0$ be the constant involved in (2.7). Define:

$$\varphi := -q - \gamma\left(1+|r|\|k\|_2^2\right) + p\frac{(1+\varepsilon)\|k\|_2^2 + (3\varepsilon-1)k_1^2}{(1+\varepsilon)\|k\|_2^2 - k_1^2}\pi^2$$
$$-p\|k\|_2^2 \max\left(0, r^2\left(\|k'\|_2^2 + (3\varepsilon-1)\pi^2\|k\|_2^2 - k'(1)\right) + r\left(2(3\varepsilon-1)\pi^2 + \frac{q}{p} - \mu\right) + \frac{3\pi^2\varepsilon(1+\varepsilon)}{(1+\varepsilon)\|k\|_2^2 - k_1^2}\right) \tag{4.21}$$

It follows from (2.7) that $\varphi > 0$. Moreover, it follows from (4.20) and the Cauchy-Schwarz inequality that the following differential inequality holds for all $t \in (0, t_{\max})$:

$$\frac{d}{dt}V(u[t]) \leq -\varphi\|u[t]\|_2^2 \tag{4.22}$$

Exploiting (4.14) and defining $\sigma := \frac{\varphi}{2c_2}$, we get from direct integration of (4.22) for all $t_0 \in (0, t_{\max})$ and $t \in [t_0, t_{\max})$:

$$V(u[t]) \leq \exp(-2\sigma(t-t_0))V(u[t_0]) \tag{4.23}$$

Estimate (2.8) with $G := \sqrt{\frac{c_2}{c_1}}$ is a direct consequence of estimates (4.23) and (4.14) and the fact that $V(u[t])$ is continuous at $t=0$ (which allows the derivation of (4.23) with $t_0 = 0$). The proof is complete. ◁

As remarked in the Introduction a CLF feedback design may not be sufficient for the derivation of stability properties in the nonlinear infinite-dimensional case. Additional analysis may be required in order to obtain bounds that guarantee the existence of the solution for all times. This crucial step is performed by the following proposition.

**Proposition 4.5:** *Let $t_{\max} \in (0,+\infty]$ and $u \in C^0([0,t_{\max}) \times [0,1]) \cap C^1((0,t_{\max}); L^2(0,1))$ with $u[t] \in H^2(0,1)$ for all $t \in [0, t_{\max})$ be a solution of (2.1), (2.2), (2.6), where $k \in C^2([0,1])$ is a function with $k(0) = 0$, $k(1) = 1$ that satisfies $k''(x) = \mu k(x)$ for $x \in [0,1]$ and for some constant $\mu \in \Re$, $p > 0$, $r \in \Re$ are constants with $1 + \|k\|_2^2 r > 0$, and $F: C^0([0,1]) \to C^0([0,1])$ is defined by (2.3) for a locally Lipschitz function $f \in C^0([0,1] \times \Re)$ that satisfies (2.4), (2.5) for certain real constants $q \in \Re$, $\gamma, \delta, B \geq 0$ and $b \geq 2$ with $B \geq \delta|r|\|k\|_b \|k\|_{\frac{b}{b-1}}$. Suppose that there exists a constant $\varepsilon \geq 0$ that satisfies (2.7). Moreover, if $q > 0$ then suppose that $B > 0$ and $b > 2$. Then the following estimate holds for all $t \in [0, t_{\max})$*

$$\|u[t]\|_\infty \leq \max\left\{\bar{K}, \|u[0]\|_\infty, G|r|\|k\|_2 \|u[0]\|_2\right\}, \tag{4.24}$$

*where $\bar{K} := (B^{-1}q)^{1/(b-2)}$ when $q > 0$ and $\bar{K} := 0$ when $q \leq 0$ and $G > 0$ is the constant involved in (2.8).*

**Proof:** We use Stampacchia's truncation method (as presented in the proof of Theorem 10.3 in [7]) as well as the fact that for every $u \in C^0([0,t_{\max}) \times [0,1]) \cap C^1((0,t_{\max}); L^2(0,1))$ and for every globally Lipschitz function $g \in C^1(\Re)$ with a globally Lipschitz derivative $g' \in C^0(\Re)$, the function $h(t) = \int_0^1 g(u(t,x))dx$ is of class $C^0([0,t_{\max})) \cap C^1((0,t_{\max}))$ with $\dot{h}(t) = \int_0^1 g'(u(t,x))u_t(t,x)dx$ for $t \in (0, t_{\max})$.

Let $t_{\max} \in (0, +\infty]$ and $u \in C^0([0,t_{\max}) \times [0,1]) \cap C^1((0,t_{\max}); L^2(0,1))$ with $u[t] \in H^2(0,1)$ for all $t \in [0, t_{\max})$ be a solution of (2.1), (2.2), (2.6). Define $M := \max\{\bar{K}, \|u[0]\|_\infty, G|r|\|k\|_2 \|u[0]\|_2\}$ and

$$T := \sup\left\{t \in [0, t_{\max}): \max_{(s,x) \in [0,t] \times [0,1]}(u(s,x)) < 1 + M\right\} \tag{4.25}$$

$$g(s) := \begin{cases} 0 & \text{if } s \leq M \\ (s-M)^3 & \text{if } s \in (M, M+1) \\ 3(s-M)-2 & \text{if } s \geq M+1 \end{cases} \tag{4.26}$$

$$h(t) := \int_0^1 g(u(t,x))dx \tag{4.27}$$

Notice that $g \in C^1(\mathfrak{R})$ as defined by (4.26) is globally Lipschitz function with a globally Lipschitz derivative $g' \in C^0(\mathfrak{R})$. Thus $h$ as defined by (4.27) is of class $C^0([0,t_{max})) \cap C^1((0,t_{max}))$ with $\dot{h}(t) = \int_0^1 g'(u(t,x))u_t(t,x)dx$ for $t \in (0,t_{max})$. By continuity of $u$ (and the maximum theorem; see page 306 in [29]), it follows from definition (4.25) that $T \in (0, t_{max}]$. Moreover, definition (4.25) implies that $u(t,x) < 1+M$ for all $(t,x) \in [0,T) \times [0,1]$. Finally, if $T < t_{max}$ then $\max_{x \in [0,1]}(u(T,x)) = 1+M$.

We show next that $T = t_{max}$. Indeed, since $\dot{h}(t) = \int_0^1 g'(u(t,x))u_t(t,x)dx$ for $t \in (0,T)$, we get from (2.1) and (2.3):

$$\dot{h}(t) = p\int_0^1 g'(u(t,x))u_{xx}(t,x)dx + \int_0^1 g'(u(t,x))f(x,u(t,x))dx \text{ for } t \in (0,T) \tag{4.28}$$

Furthermore, notice that Proposition 4.4 (estimate (2.8)) in conjunction with definition $M := \max\{\bar{K}, \|u[0]\|_\infty, G|r|\|k\|_2 \|u[0]\|_2\}$, (2.6) and the Cauchy-Schwarz inequality implies that $u(t,1) \leq M$ for $t \in (0,T)$. Therefore, we get (using (2.6) and definition (4.26)) that $g'(u(t,0)) = g'(u(t,1)) = 0$ for $t \in (0,T)$. The fact $u(t,x) < 1+M$ for all $(t,x) \in [0,T) \times [0,1]$ in conjunction with the fact that $g$ is of class $C^2$ on $(-\infty, M+1)$ allows us to use integration by parts (notice that $u[t] \in H^2(0,1)$ for all $t \in [0, t_{max})$) and write (4.28) in the following way:

$$\dot{h}(t) = -p\int_0^1 g''(u(t,x))(u_x(t,x))^2 dx + \int_0^1 g'(u(t,x))f(x,u(t,x))dx \text{ for } t \in (0,T) \tag{4.29}$$

The fact that $u(t,x) < 1+M$ for all $(t,x) \in [0,T) \times [0,1]$ in conjunction with the fact that $g''(s) \geq 0$ for all $s \in (-\infty, M+1)$ and (4.29) gives:

$$\dot{h}(t) \leq \int_0^1 g'(u(t,x))f(x,u(t,x))dx \text{ for } t \in (0,T) \tag{4.30}$$

Moreover, inequality (4.30) in conjunction with the facts that:
- $u(t,x) < 1+M$ for all $(t,x) \in [0,T) \times [0,1]$,
- $g'(u(t,x)) = 0$ for all $(t,x) \in [0,T) \times [0,1]$ with $u(t,x) \leq M$ (a consequence of definition (4.26)),
- $g'(u(t,x)) \geq 0$ for all $(t,x) \in [0,T) \times [0,1]$ with $M < u(t,x) < 1+M$ (a consequence of definition (4.26)),
- $f(x,u(t,x)) \leq 0$ for all $(t,x) \in [0,T) \times [0,1]$ with $M < u(t,x) < 1+M$ (a consequence of definitions $M := \max\{\bar{K}, \|u[0]\|_\infty, G|r|\|k\|_2 \|u[0]\|_2\}$, $\bar{K} := (B^{-1}q)^{1/(b-2)}$ when $q > 0$ and $\bar{K} := 0$ when $q \leq 0$ and inequality (2.4)),

allows us to conclude that

$$\dot{h}(t) \leq 0 \text{ for } t \in (0,T) \tag{4.31}$$

It follows from (4.31) and continuity of $h$ that $h(t) \leq h(0)$ for $t \in [0,T)$. Since $h(0) = 0$ (recall (4.26), (4.27) and definition $M := \max\{\bar{K}, \|u[0]\|_\infty, G|r|\|k\|_2 \|u[0]\|_2\}$) we get that $h(t) \leq 0$ for $t \in [0,T)$. Continuity of $u$ and definitions (4.26), (4.27) imply that $u(t,x) \leq M$ for all $(t,x) \in [0,T) \times [0,1]$.

As remarked earlier, if $T < t_{\max}$ then $\max_{x \in [0,1]}(u(T,x)) = 1 + M$; a contradiction with continuity of $u$ and the fact that $u(t,x) \le M$ for all $(t,x) \in [0,T) \times [0,1]$. Therefore $T = t_{\max}$ and $u(t,x) \le M$ for all $(t,x) \in [0,t_{\max}) \times [0,1]$.

A similar analysis is used in order to show that $u(t,x) \ge -M$ for all $(t,x) \in [0,t_{\max}) \times [0,1]$. Thus (4.24) holds. The proof is complete. ◁

**Remark 4.6:** The proof of Proposition 4.5, which uses Proposition 4.4, shows another one important technical issue that would arise in the case of a nonlinear boundary feedback law: since Proposition 4.4 provides a bound for the $L^2$ spatial norm of the state, the boundary control action given by a nonlinear feedback law would have to be bounded by expressions involving only the $L^2$ spatial norm of the state.

Combining Theorem 3.1 and Proposition 4.5, we obtain the following corollary.

**Corollary 4.7:** *Consider the closed-loop system (2.1), (2.2), (2.6), where $k \in C^2([0,1])$ is a function with $k(0) = 0$, $k(1) = 1$ that satisfies $k''(x) = \mu k(x)$ for $x \in [0,1]$ and for some constant $\mu \in \Re$, $p > 0$, $r \in \Re$ are constants with $1 + \|k\|_2^2 r > 0$, and $F : C^0([0,1]) \to C^0([0,1])$ is defined by (2.3) for a locally Lipschitz function $f \in C^0([0,1] \times \Re)$ that satisfies (2.4), (2.5) for certain real constants $q \in \Re$, $\gamma, \delta, B \ge 0$ and $b \ge 2$ with $B \ge \delta |r| \|k\|_b \|k\|_{\frac{b}{b-1}}$. Suppose that there exists a constant $\varepsilon \ge 0$ that satisfies (2.7). Moreover, if $q > 0$ then suppose that $B > 0$ and $b > 2$. Then for every $u_0 \in H^2(0,1)$ with $u_0(0) = 0$, $u_0(1) = -r \langle k, u_0 \rangle$, there exists a unique mapping $u \in C^0(\Re_+ \times [0,1]) \cap C^1((0,+\infty); L^2(0,1))$ with $u[0] = u_0$, $u[t] \in H^2(0,1)$ for all $t \ge 0$ for which equation (2.1) holds for all $t > 0$ and equations (2.2), (2.6) hold for all $t \ge 0$. Moreover, there exist constants $G, \sigma > 0$ such that estimates (2.8), (4.24) hold for all $t \ge 0$. Finally, the mapping $t \to \|u_x[t] + rku[t]\|_2^2$ is $C^1$ on $\Re_+$ and the following equation holds*

$$\frac{d}{dt} \|u_x[t] + rku[t]\|_2^2 = -2p \|u_{xx}[t] + rku_x[t] + rk'u[t]\|_2^2 - 2 \langle KF(u[t]) + pGu[t], u_{xx}[t] + rku_x[t] + rk'u[t] \rangle \quad (4.32)$$

*where $K, G : L^2(0,1) \to L^2(0,1)$ are the continuous linear operators defined by the following equations for all $u \in L^2(0,1)$, $x \in [0,1]$:*

$$(Ku)(x) = u(x) + r \int_0^x k(s)u(s)ds \quad (4.33)$$

$$(Gu)(x) = -2rk'(x)u(x) + r \int_0^x k''(s)u(s)ds \quad (4.34)$$

We are now ready to give the proof of the main result of the paper.

**Proof of Theorem 2.1:** Using (2.3), (2.5), the triangle inequality and the fact that $b \ge 2$ we get

$$\|KF(u) + pGu\|_2 \le \left( p\|G\|_2 + \|K\|_2 \left( \gamma + |q| + \delta \|u\|_\infty^{b-2} \right) \right) \|u\|_2 \text{ for all } w \in C^0([0,1]) \quad (4.35)$$

where $\|G\|_2 = \sup\{ \|Gv\|_2 : v \in L^2(0,1), \|v\|_2 = 1 \}$, $\|K\|_2 = \sup\{ \|Kv\|_2 : v \in L^2(0,1), \|v\|_2 = 1 \}$.

Let $u_0 \in H^2(0,1)$ with $u_0(0) = 0$, $u_0(1) = -r \langle k, u_0 \rangle$, be given (arbitrary) and consider the unique mapping $u \in C^0(\Re_+ \times [0,1]) \cap C^1((0,+\infty); L^2(0,1))$ with $u[0] = u_0$, $u[t] \in H^2(0,1)$ for all $t \ge 0$ for which equation (2.1) holds for all $t > 0$ and equations (2.2), (2.6) hold for all $t \ge 0$ (whose existence and uniqueness is guaranteed by Corollary 4.7).

Using (4.32), (4.35) and the Cauchy-Schwarz inequality, we get for $t \ge 0$

$$\frac{d}{dt}\|u_x[t]+rku[t]\|_2^2 \leq -2p\|u_{xx}[t]+rku_x[t]+rk'u[t]\|_2^2 + 2g(\|u[t]\|_\infty)\|u_{xx}[t]+rku_x[t]+rk'u[t]\|_2 \|u[t]\|_2 \quad (4.36)$$

where $g(s) := p\|G\|_2 + \|K\|_2(\gamma+|q|+\delta s^{b-2})$ for $s \geq 0$ is a non-decreasing function. Using the fact that

$$2g(\|u[t]\|_\infty)\|u_{xx}[t]+rku_x[t]+rk'u[t]\|_2 \|u[t]\|_2 \leq 2p\|u_{xx}[t]+rku_x[t]+rk'u[t]\|_2^2 + \bar{g}(\|u[t]\|_\infty)\|u[t]\|_2^2,$$

where $\bar{g}(s) := \frac{1}{2p}(g(s))^2$ for $s \geq 0$ is a non-decreasing function, we get from (4.36):

$$\frac{d}{dt}\|u_x[t]+rku[t]\|_2^2 \leq \bar{g}(\|u[t]\|_\infty)\|u[t]\|_2^2, \text{ for } t \geq 0 \quad (4.37)$$

Using (2.8) in conjunction with (4.37) and (4.24) we get:

$$\frac{d}{dt}\|u_x[t]+rku[t]\|_2^2 \leq G^2 \bar{g}(M) \exp(-2\sigma t)\|u_0\|_2^2, \text{ for } t \geq 0 \quad (4.38)$$

where $M := \max\{\bar{K}, \|u[0]\|_\infty, G|r|\|k\|_2 \|u[0]\|_2\}$, $\bar{K} := (B^{-1}q)^{1/(b-2)}$ when $q>0$ and $\bar{K} := 0$ when $q \leq 0$ and $G, \sigma > 0$ are the constants involved in (2.8). The differential inequality (4.38) directly implies

$$\|u_x[t]+rku[t]\|_2^2 \leq \|u_0'+rku_0\|_2^2 + \frac{1}{2\sigma}G^2 \bar{g}(M)\|u_0\|_2^2, \text{ for } t \geq 0$$

or

$$\|u_x[t]+rku[t]\|_2 \leq \|u_0'+rku_0\|_2 + \tilde{g}(M)\|u_0\|_2, \text{ for } t \geq 0 \quad (4.39)$$

where $\tilde{g}(s) := G\sqrt{\frac{\bar{g}(s)}{2\sigma}}$ for $s \geq 0$ is a non-decreasing function. Using the facts that $\|u_x[t]\|_2 \leq \|u_x[t]+rku[t]\|_2 + |r|\|k\|_\infty \|u[t]\|_2$, $\|u_0'+rku_0\|_2 \leq \|u_0'\|_2 + |r|\|k\|_\infty \|u_0\|_2$ (both consequences of the triangle inequality), we obtain from (4.39) in conjunction with (2.8):

$$\|u_x[t]\|_2 \leq \|u_0'\|_2 + \psi(M)\|u_0\|_2, \text{ for } t \geq 0 \quad (4.40)$$

where $\psi(s) := \tilde{g}(s) + (1+G)|r|\|k\|_\infty$ for $s \geq 0$ is a non-decreasing function. Estimate (2.10) is a direct consequence of estimates (4.24) and (4.40). Finally, estimate (2.9) is a direct consequence of estimates (2.8), (2.10) and Agmon's inequality

$$\|u\|_\infty \leq \sqrt{2}\sqrt{\|u\|_2 \|u'\|_2} \text{ , for all } u \in H^1(0,1) \text{ with } u(0)=0$$

The proof is complete. ◁

## 5. Concluding Remarks

The present paper showed what difficulties may be encountered in the extension of well-known feedback design methodologies for finite-dimensional systems to PDEs. More specifically, it was shown that in the nonlinear infinite-dimensional case a CLF feedback design may not be sufficient for establishing existence of solutions for all times (global solutions) and consequently may not allow a valid (not merely formal) derivation of stability properties. Additional analysis may be required in order to obtain bounds that guarantee the existence of the solution for all times and the pointwise convergence of the solution to the desired equilibrium point (important in practice).

It should be clear that this paper is only a first step towards the design of global boundary feedback stabilizers for 1-D nonlinear parabolic PDEs. Additional steps are needed and will be the topic of future research. For example, the stabilization of unstable PDEs with more than one unstable mode will require the development of novel mathematical results which will allow the construction of appropriate CLFs.